\documentclass[12pt,twoside]{article}

\usepackage[english]{babel}
\usepackage{umj_1}           

\usepackage{graphicx, srcltx}   
\usepackage{cite}               

\usepackage{amsmath}            
\usepackage{amsfonts,amssymb}   
\usepackage{srcltx}             

\newcommand{\Leftclt}{Roman A. Veprintsev}
\newcommand{\Rightclt}{On Paley-type and Hausdorff\,--\,Young\,--\,Paley-type inequalities}
\usepackage{authblk}

\DeclareMathOperator*{\esssup}{ess\,sup}

\begin{document}

\begin{Titul}
{\large \bf ON PALEY-TYPE\\ AND HAUSDORFF\,--\,YOUNG\,--\,PALEY-TYPE INEQUALITIES\\[0.2em] FOR JACOBI EXPANSIONS }\\[3ex]
{{\bf Roman~A.~Veprintsev} \\[5ex]}
\end{Titul}

\begin{Anot}
{\bf Abstract.} We obtain Paley-type and Hausdorff\,--\,Young\,--\,Paley-type inequalities for Jacobi expansions.


\end{Anot}


\section{Introduction and preliminaries}

In this section, we introduce the basic notation and review some results on Jacobi polynomials. Here, we also formulate the aim of the publication.

Let $\mathbb{N}_0$ denote the set of non-negative integers. Suppose that $\alpha,\,\beta>-1$. The Jacobi polynomials, denoted by $P_n^{(\alpha,\beta)}(\cdot)$, $n\in\mathbb{N}_0$, are orthogonal with respect to the Jacobi weight function $w_{\alpha,\beta}(t)=(1-t)^\alpha(1+t)^\beta$ on $[-1,1]$, namely,
\begin{equation*}
\int\nolimits_{-1}^1 P_n^{(\alpha,\beta)}(t)\, P_m^{(\alpha,\beta)}(t)\, w_{\alpha,\beta}(t)\,dt=\begin{cases}\dfrac{2^{\alpha+\beta+1}\Gamma(n+\alpha+1)\Gamma(n+\beta+1)}{(2n+\alpha+\beta+1)\Gamma(n+1)\Gamma(n+\alpha+\beta+1)},&n=m,\\
0,&n\not=m.
\end{cases}
\end{equation*}
Here, as usual, $\Gamma$ is the gamma function.

We use the following asymptotic notation: $f(n)\asymp g(n)$, $n\to\infty$, means that there exist positive constants $C_1$, $C_2$, and a positive integer $n_0$ such that $0\leq C_1 g(n)\leq f(n)\leq C_2 g(n)$ for all $n\geq n_0$.
To simplify the writing, we will omit ``$n\to\infty$'' in the notation.

It follows directly from Stirling's asymptotic formula that
$
\frac{\Gamma(n+\lambda)}{\Gamma(n+\mu)}\asymp n^{\lambda-\mu}
$
for arbitrary real numbers $\lambda$ and $\mu$.
Thus, we have
\begin{equation}\label{asymptotic_for_L_2-norm_of_Jacobi_polynomials}
h_n^{(\alpha,\beta)}=\Bigl(\int\nolimits_{-1}^1 \bigl(P_n^{(\alpha,\beta)}(t)\bigr)^2\, w_{\alpha,\beta}(t)\,dt\Bigr)^{1/2}\asymp n^{-1/2}.
\end{equation}

Define the uniform norm of a continuous function $f$ on $[-1,1]$ by
\begin{equation*}
\|f\|_{\infty}=\max\limits_{-1\leq t\leq 1} |f(t)|.
\end{equation*}
The maximum of two real numbers $x$ and $y$ is denoted by $\max(x,y)$.

It is known \cite[Theorem~7.32.1]{szego_book_orthogonal_polynomials_1975} that
\begin{equation}\label{asymptotic_for_maximum_of_absolute_value_of_Jacobi_polynomials}
\bigl\|P_n^{(\alpha,\beta)}\bigr\|_\infty\asymp n^{\sigma(\alpha,\beta)-\frac{1}{2}},\qquad \text{where}\quad \sigma(\alpha,\beta)=\max\Bigl(0,\frac{1}{2}+\max(\alpha,\beta)\Bigr).
\end{equation}

Let $\bigl\{\widetilde{P}_n^{(\alpha,\beta)}\bigr\}_{n=0}^\infty$ be the sequence of orthonormal Jacobi polynomials (with respect to the Jacobi weight function $w_{\alpha,\beta}$), that is,  $\widetilde{P}_n^{(\alpha,\beta)}\hm{=}\bigl(h_n^{(\alpha,\beta)}\bigr)^{-1} P_n^{(\alpha,\beta)}$, $n\in\mathbb{N}_0$.

It follows from \eqref{asymptotic_for_L_2-norm_of_Jacobi_polynomials} and \eqref{asymptotic_for_maximum_of_absolute_value_of_Jacobi_polynomials} that
\begin{equation}\label{asymptotic_for_maximum_of_absolute_value_of_orthonormal_Jacobi_polynomials}
\bigl\|\widetilde{P}_n^{(\alpha,\beta)}\bigr\|_\infty\asymp n^{\sigma(\alpha,\beta)}.
\end{equation}

Given $1\leq p\leq\infty$, we denote by $L_p(w_{\alpha,\beta})$ the space of complex-valued Lebesgue measurable functions $f$ on $[-1,1]$ with finite norm
\begin{equation*}
\begin{array}{lr}
\|f\|_{L_p(w_{\alpha,\beta})}=\Bigl(\int\nolimits_{-1}^1 |f(t)|^p\,w_{\alpha,\beta}(t)\,dt\Bigr)^{1/p},&\quad 1\leq p<\infty,\\[1.0em]
\|f\|_{L_\infty}=\esssup\limits_{x\in[-1,1]} |f(x)|,& p=\infty.
\end{array}
\end{equation*}
For a function $f\in L_p(w_{\alpha,\beta})$, $1\leq p\leq\infty$, the Jacobi expansion is given by
\begin{equation*}
f(t)\sim\sum\limits_{n=0}^\infty \hat{f}_n \widetilde{P}_n^{(\alpha,\beta)}(t),\qquad \text{where}\quad\hat{f}_n=\int\nolimits_{-1}^1 f(t)\, \widetilde{P}_n^{(\alpha,\beta)}(t)\,w_{\alpha,\beta}(t)\,dt.
\end{equation*}

For $1\leq p\leq\infty$, we denote by $p'$ the conjugate exponent to $p$, that is, $\frac{1}{p}+\frac{1}{p'}=1$.

Our aim is to establish Paley-type and Hausdorff\,--\,Young\,--\,Paley-type inequalities for Jacobi expansions in Sections \ref{section_for_Paley-type_inequality} and \ref{section_for_Hausdorff-Young-Paley-type_inequality}, respectively.

\section{Paley-type inequalities for Jacobi expansions}\label{section_for_Paley-type_inequality}

We will omit the proof of the following Paley-type inequalities for Jacobi expansions because its proof is analogous to that of Theorem 2 given in \cite[Section~2]{veprintsev_preprint_Paley_inequalities_2016}.

\begin{teoen}\label{main_theorem_Paley-type_inequality_for_Jacobi_expansions}
$(a)$ If $1<p\leq 2$, $f\in L_{p}(w_{\alpha,\beta})$, $\omega$ is a positive function on $\mathbb{N}_0$ such that
\begin{equation}\label{condition_in_Paley_inequality}
M_{\omega}=\sup\limits_{t>0} \Bigl\{t\Bigl(\sum\limits_{\omega(n)\geq t} (n+1)^{2\sigma(\alpha,\beta)}\Bigr)\Bigr\}<\infty,
\end{equation}
then
\begin{equation}\label{first_part_of_Paley_inequality}
\Bigl\{\sum\limits_{n=0}^\infty \Bigl((n+1)^{\left(\frac{1}{p}-\frac{1}{p'}\right)\sigma(\alpha,\beta)}\,(\omega(n))^{\frac{1}{p}-\frac{1}{p'}}|\hat{f}_n|\Bigr)^p\Bigr\}^{1/p}\leq A_p M_{\omega}^{\frac{1}{p}-\frac{1}{p'}} \|f\|_{L_p(w_{\alpha,\beta})}.
\end{equation}

$(b)$ If $2\leq q<\infty$, $\omega$ is a positive function on $\mathbb{N}_0$ satisfying \eqref{condition_in_Paley_inequality} and $\phi$ is a non-negative function on $\mathbb{N}_0$ such that
\begin{equation*}\label{assumption_for_second_part_of_Paley_inequality}
\sum\limits_{n=0}^\infty \Bigl((n+1)^{\left(\frac{1}{q}-\frac{1}{q'}\right)\sigma(\alpha,\beta)}(\omega(n))^{\frac{1}{q}-\frac{1}{q'}}|\phi(n)|\Bigr)^{q}<\infty,
\end{equation*}
then the algebraic polynomials
\begin{equation*}
\Phi_N(t)=\sum\limits_{n=0}^N \phi(n)\widetilde{P}_n^{(\alpha,\beta)}(t)
\end{equation*}
converge in $L_q(w_{\alpha,\beta})$ to a function $f$ satisfying $\hat{f}_n=\phi(n)$, $n\in\mathbb{N}_0$, and
\begin{equation*}\label{second_part_of_Paley_inequality}
\|f\|_{L_q(w_{\alpha,\beta})}\leq A_{q'} M_{\omega}^{\frac{1}{q'}-\frac{1}{q}} \Bigl\{\sum\limits_{n=0}^\infty \Bigl((n+1)^{\left(\frac{1}{q}-\frac{1}{q'}\right)\sigma(\alpha,\beta)}(\omega(n))^{\frac{1}{q}-\frac{1}{q'}}|\phi(n)|\Bigr)^{q}\Bigr\}^{1/q}.
\end{equation*}
\end{teoen}

\section{Hausdorff\,--\,Young\,--\,Paley-type inequalities\\ for Jacobi expansions}\label{section_for_Hausdorff-Young-Paley-type_inequality}

In \cite[Theorem~3.1]{ditzian_article_norm_and_smoothness_2015}, Z.~Ditzian proved the following Hausdorff\,--\,Young-type inequalities for Jacobi expansions.

\begin{teoen}\label{Hausdorff-Young-Paley_inequality_for_generalized_Gegenbauer_expansions}
$(a)$ If $1<p\leq2$ and $f\in L_{p}(w_{\alpha,\beta})$, then
\begin{equation*}\label{first_part_of_Hausdorff-Young_inequality}
\Bigl\{\sum\limits_{n=0}^\infty\Bigl((n+1)^{\left(\frac{1}{p'}-\frac{1}{p}\right)\sigma(\alpha,\beta)}|\hat{f}_n|\Bigr)^{p'}\Bigr\}^{1/p'}\leq B_p\,\|f\|_{L_p(w_{\alpha,\beta})}.
\end{equation*}

$(b)$ If $2\leq q<\infty$ and $\phi$ is a function on non-negative integers satisfying
\begin{equation*}\label{assumption_for_second_part_of_Hausdorff-Young_inequality}
\sum\limits_{n=0}^\infty\Bigl((n+1)^{\left(\frac{1}{q'}-\frac{1}{q}\right)\sigma(\alpha,\beta)}|\phi(n)|\Bigl)^{q'}<\infty,
\end{equation*}
then the algebraic polynomials
\begin{equation*}
\Phi_N(t)=\sum\limits_{n=0}^N \phi(n)\,\widetilde{P}_n^{(\alpha,\beta)}(t)
\end{equation*}
converge in $L_q(w_{\alpha,\beta})$ to a function $f$ satisfying $\hat{f}_n=\phi(n)$, $n\in\mathbb{N}_0$, and
\begin{equation*}\label{second_part_of_Hausdorff-Young_inequality}
\|f\|_{L_q(w_{\alpha,\beta})}\leq B_{q'} \Bigl\{\sum\limits_{n=0}^\infty\Bigl((n+1)^{\left(\frac{1}{q'}-\frac{1}{q}\right)\sigma(\alpha,\beta)}|\phi(n)|\Bigl)^{q'}\Bigr\}^{1/q'}.
\end{equation*}
\end{teoen}

Theorem \ref{unification_of_different_types_of_inequalities} contains the Paley-type and the Hausdorff\,--\,Young-type inequalities for Jacobi expansions. The proof of this theorem is analogous to the proof of Theorem 5 in \cite[Section~3]{veprintsev_preprint_Paley_inequalities_2016} and we shall omit it.

\begin{teoen}\label{unification_of_different_types_of_inequalities}
$(a)$ If $1<p\leq 2$, $f\in L_p(w_{\alpha,\beta})$, $\omega$ is a positive function on $\mathbb{N}_0$ satisfying the condition \eqref{condition_in_Paley_inequality}, and $p\leq s\leq p'$, then
\begin{equation*}\label{first_part_of_unification}
\Bigl\{\sum\limits_{n=0}^\infty\Bigl((n+1)^{\left(\frac{2}{s}-1\right)\sigma(\alpha,\beta)}(\omega(n))^{\frac{1}{s}-\frac{1}{p'}}|\hat{f}_n|\Bigr)^s\Bigr\}^{1/s}\leq C_p(s)\,M_\omega^{\frac{1}{s}-\frac{1}{p'}} \,\|f\|_{L_p(w_{\alpha,\beta})}.
\end{equation*}

$(b)$ If $2\leq q<\infty$, $q'\leq r\leq q$, $\omega$ is a positive function on $\mathbb{N}_0$ satisfying the condition \eqref{condition_in_Paley_inequality} and $\phi$ is a non-negative function on $\mathbb{N}_0$ such that
\begin{equation*}
\sum\limits_{n=0}^\infty \Bigl((n+1)^{\left(1-\frac{2}{r}\right)\sigma(\alpha,\beta)}(\omega(n))^{\frac{1}{q}-\frac{1}{r}}|\phi(n)|\Bigr)^{r'}<\infty,
\end{equation*}
then the algebraic polynomials
\begin{equation*}
\Phi_N(t)=\sum\limits_{n=0}^N \phi(n)\widetilde{P}_n^{(\alpha,\beta)}(t)
\end{equation*}
converge in $L_q(w_{\alpha,\beta})$ to a function $f$ satisfying $\hat{f}_n=\phi(n)$, $n\in\mathbb{N}_0$, and
\begin{equation*}
\|f\|_{L_q(w_{\alpha,\beta})}\leq C_{q'}(r) M_{\omega}^{\frac{1}{r}-\frac{1}{q}} \Bigl\{\sum\limits_{n=0}^\infty \Bigl((n+1)^{\left(1-\frac{2}{r}\right)\sigma(\alpha,\beta)}(\omega(n))^{\frac{1}{q}-\frac{1}{r}}|\phi(n)|\Bigr)^{r'}\Bigr\}^{1/r'}.
\end{equation*}
\end{teoen}

\section{Remarks}

In the following theorem we show that the inequality \eqref{first_part_of_Paley_inequality} does not hold for $p=1$ in general.

\begin{teoen}
Under the conditions and notations of {\rm Theorem \ref{main_theorem_Paley-type_inequality_for_Jacobi_expansions}}, \eqref{first_part_of_Paley_inequality} is not satisfied for $p=1$ provided that $\max(\alpha,\beta)\geq-\frac{1}{2}$ and
\begin{equation*}
\sum\limits_{n=0}^\infty \omega(n)(n+1)^{2\sigma(\alpha,\beta)}=+\infty.
\end{equation*}
\end{teoen}

\proofen We assume without loss of generality that $\max(\alpha,\beta)=\alpha$. Then we have $\sigma(\alpha,\beta)=\frac{1}{2}+\alpha$. We define the sequence of polynomials $g_N$ by \[g_N(t)=\sum\limits_{n=0}^N \omega(n) (n+1)^{\sigma(\alpha,\beta)} \widetilde{P}_n^{(\alpha,\beta)}(t).\]
As $\bigl\|\widetilde{P}_n^{(\alpha,\beta)}\bigr\|_\infty=\widetilde{P}_n^{(\alpha,\beta)}(1)\asymp (n+1)^{\sigma(\alpha,\beta)}$ (see \cite[Theorem~7.32.1, (4.1.1)]{szego_book_orthogonal_polynomials_1975} and \eqref{asymptotic_for_maximum_of_absolute_value_of_orthonormal_Jacobi_polynomials}), we obtain $\|g_N\|_\infty=\sum\limits_{n=0}^N\omega(n)(n+1)^{\sigma(\alpha,\beta)} \widetilde{P}_n^{(\alpha,\beta)}(1)\asymp\sum\limits_{n=0}^N \omega(n)(n+1)^{2\sigma(\alpha,\beta)}$. Thus,
\begin{equation*}
\|g_N\|_\infty\to\infty\quad\text{as}\quad N\to\infty.
\end{equation*}

We now show that assuming \eqref{first_part_of_Paley_inequality} for $p=1$ would lead to a contradiction with what we already proved. Because of $\|g_N\|_\infty=\sup\,\bigl\{\int\nolimits_{-1}^1 g_N(t) f(t)\,w_{\alpha,\beta}(t)\,dt\colon\, \|f\|_{L_1(w_{\alpha,\beta})}=1\bigr\}$, we have
\begin{equation*}
\begin{split}
\|g_N\|_\infty&=\sup\, \Bigl\{\int\nolimits_{-1}^1 g_N(t) f(t)\,w_{\alpha,\beta}(t)\,dt\colon\, \|f\|_{L_1(w_{\alpha,\beta})}=1\Bigr\}=\\
&=\sup\, \Bigl\{\sum\limits_{n=0}^N \omega(n) (n+1)^{\sigma(\alpha,\beta)} \hat{f}_n\colon\, \|f\|_{L_1(w_{\alpha,\beta})}=1\Bigr\}\leq\\
&\leq \sup\, \Bigl\{\sum\limits_{n=0}^N \omega(n) (n+1)^{\sigma(\alpha,\beta)} |\hat{f}_n|\colon\, \|f\|_{L_1(w_{\alpha,\beta})}=1\Bigr\}\leq A_1 M_\omega.
\end{split}
\end{equation*}
This is a contradiction, which proves that our assumption is wrong.
\hfill$\square$

\section{Conclusion}

As an application of the results above, we are going to obtain sufficient condition for a special sequence of positive real numbers to be a Fourier multiplier.

\section*{Acknowledgements}

The author would like to thank Michael Ruzhansky for the idea of this work (see the e-print \cite{akylzhanov_nursultanov_ruzhansky_article_inequalities_2015}). This work was done thanks to the remarkable paper \cite{ditzian_article_norm_and_smoothness_2015} (see especially Section 4 there) of Z.~Ditzian.

\begin{Biblioen}

\bibitem{akylzhanov_nursultanov_ruzhansky_article_inequalities_2015}R.~Akylzhanov, E.~Nursultanov, and M.~Ruzhansky, Hardy\,--\,Littlewood, Hausdorff\,--\,Young\,--\,Paley inequalities, and $L^p$-$L^q$ Fourier multipliers on compact homogeneous manifolds, arXiv preprint 1504.07043 (2015).





\bibitem{ditzian_article_norm_and_smoothness_2015}Z.~Ditzian, Norm and smoothness of a function related to the coefficients of its expansion, \textit{J. Approx. Theory} \textbf{196} (2015), 101--110.



\bibitem{stein_article_interpolation_1956}E.~Stein, Interpolation of linear operators, \textit{Trans. Amer. Math. Soc.} \textbf{83} (1956), 482--492.


\bibitem{szego_book_orthogonal_polynomials_1975}G.~Szeg\"{o}, \textit{Orthogonal polynomials}, 4th ed., American Mathematical Society Colloquium Publications \textbf{23}, American Mathematical Society, Providence, Rhode Island, 1975.



\bibitem{veprintsev_preprint_Paley_inequalities_2016}R.\,A.~Veprintsev, On Paley-type and Hausdorff\,--\,Young\,--\,Paley-type inequalities for generalized Gegenbauer expansions, arXiv preprint 1603.02072 (2016).

\end{Biblioen}

\noindent \textsc{Independent researcher, Uzlovaya, Russia}

\noindent \textit{E-mail address}: \textbf{veprintsevroma@gmail.com}

\end{document}